\documentclass[a4paper,11pt]{amsart}

\usepackage{amsmath,amssymb,amsthm}

\numberwithin{equation}{section} \theoremstyle{plain}
\newtheorem{thm}{Theorem}[section]

\newtheorem{lem}[thm]{Lemma}

\newtheorem{rem}[thm]{Remark}

\newtheorem{ques}[thm]{Question}

\newtheorem*{Data}{Data availability statement}
\newtheorem*{acknow}{Acknowledgments}
   
 \makeatletter

\def\<{\langle}
\def\>{\rangle}
\def\({\left(}
\def\){\right)}
\def\[{\left[}
\def\]{\right]}

\makeatother
\title {Cohomology vanishing theorems for free boundary submanifolds}
\author[N. Chen]{Niang Chen}
\address{School of Mathematical Sciences, Laboratory of Mathematics and Complex Systems, Beijing Normal University, Beijing 100875, P.R. CHINA.}
\email{chenniangbnu@mail.bnu.edu.cn}

\author[J.Q. Ge]{Jianquan Ge$^{*}$}
\address{School of Mathematical Sciences, Laboratory of Mathematics and Complex Systems, Beijing Normal University, Beijing 100875, P.R. CHINA.}
\email{jqge@bnu.edu.cn}

\date{}
\thanks {$^{*}$ the corresponding author.}
\thanks{The second author is partially supported by Beijing Natural Science Foundation (No. Z190003) and NSFC (No. 12171037).}

\begin{document}
\maketitle

\begin{abstract}
In this paper, via a new Hardy type inequality, we establish some cohomology vanishing theorems for free boundary compact submanifolds $M^n$ with $n\geq2$ immersed in the Euclidean unit ball $\mathbb{B}^{n+k}$ under one of the pinching conditions $|\Phi|^2\leq C$, $|A|^2\leq \widetilde{C}$, or $|\Phi|\leq R(p,|H|)$, where $A$ $(\Phi)$ is the (traceless) second fundamental form, $H$ is the mean curvature, $C,\widetilde{C}$ are positive constants and $R(p,|H|)$ is a positive function. In particular, we remove the condition on the flatness of the normal bundle, solving the first question, and partially answer the second question on optimal pinching constants proposed by Cavalcante, Mendes and Vit\'{o}rio.
\end{abstract}
\section{Introduction}

Relations between geometry and topology are always fascinating to geometers. Initiated by Simons's pinching rigidity theorems (cf. \cite{Simons68, CDK1970, Lawson1969}, etc.), extensive research has been done under pinching conditions for the second fundamental form $A$ of submanifolds. In 1973, Lawson and Simons \cite{Lawson-Simons} proved that if $M^n$ is a closed submanifold immersed in the unit sphere $\mathbb{S}^{n+k}$ satisfying $|A|^2<\min\{p(n-p),2\sqrt{p(n-p)}\}$ where $p<n$ is a positive integer, then the homology groups $H_p(M;G)=H_{n-p}(M;G)=0$ for any finitely generated Abelian group. Namely, bounds on the length of the second fundamental form for submanifolds of spheres imply the vanishing of homology groups. Later, plenty of generalizations with bounds involving mean curvature and various applications to differentiable sphere theorems or eigenvalue's estimates were often obtained (cf. \cite{Xin84, Vlachos, Savo 2014, XG2013, CS2019}, etc.).

It is often observed that free boundary submanifolds in unit balls have similar properties as closed submanifolds in unit spheres (cf. \cite{FL2014, AN20, CMV2019}, etc.).
 A submanifold $M^n$ with nonempty boundary $\partial M $ immersed in the unit ball $\mathbb{B}^{n+k}$ $(k\geq1)$ is called \textit{free boundary}, if  $M^n \cap \partial \mathbb{B}^{n+k} = \partial M$ and $M^n$ intersects $\partial \mathbb{B}^{n+k}=\mathbb{S}^{n+k-1}$ orthogonally along $\partial M$.
Analogous to the Lawson-Simons homology vanishing theorem for closed submanifolds of unit spheres, Cavalcante, Mendes and Vit\'{o}rio \cite{CMV2019} obtained the following cohomology vanishing theorems for compact free boundary submanifolds in unit balls, where the traceless second fundamental form $\Phi:=A-H\mathrm{Id}$ is taken in place of the second fundamental form $A$, $H:=\frac{1}{n}\mathrm{Tr}A$ is the mean curvature vector field, and $|\Phi|^2=|A|^2-n|H|^2$.

\begin{thm}\cite{CMV2019}\label{thmrefer1}
Let $M^n $ be a compact oriented submanifold immersed in $\mathbb B^{n+k}$, with $n \ge 3$, which is free boundary and has flat normal bundle. If $|\Phi|^2< \cfrac{np}{n-p}$, for some positive integer $p \leq \lfloor \frac{n}{2} \rfloor$, then the $p$th and the $(n-p)$th cohomology groups of  $M^n $ with real coefficients vanish, that is, $H^p(M; \mathbb{R}) = H^{ n-p}(M; \mathbb{R}) = 0$. In particular, if $|\Phi|^2< \frac{n}{n-1}$, then all cohomology groups $H^q(M; \mathbb{R})$, with $q = 1, \dots , n-1$, vanish and $M$ has only one boundary component.
\end{thm}
\begin{thm}\cite{CMV2019}\label{thmrefer2}
Let $M^n $ be a compact oriented submanifold minimally immersed in $\mathbb B^{n+k}$, with $n \ge 3$, which is free boundary and has flat normal bundle. Given a positive integer $p \leq \lfloor \frac{n}{2} \rfloor$, we have the following assertions:
\begin{itemize}
\item[$(1)$] If $|A|^2<\frac{n^2}{2(n-p)}$, then $H^p(M; \mathbb{R})$ vanishes. If, additionally, $p=\lfloor \frac{n}{2} \rfloor$, then $H^{n-p}(M; \mathbb{R})$ also vanishes.
\item[$(2)$] If $|A|^2\leq\frac{(n-p+1)n^3}{4p(n-p)^2}$ and $1\leq p\leq \lfloor \frac{n}{2} \rfloor-1$, then $H^{n-p}(M; \mathbb{R})$ vanishes.
\end{itemize}
In particular, if $|A|^2<\frac{n^2}{2(n-1)}$, then all cohomology groups $H^q(M; \mathbb{R})$, with $q = 1, \dots , n-1$, vanish and $M$ has only one boundary component.
\end{thm}
These theorems lead them to the following
\begin{ques}\cite{CMV2019}\label{questions}
Do Theorems $\ref{thmrefer1}$ and $\ref{thmrefer2}$ hold without the condition on the flatness of the normal bundle? What are the best constants in such theorems?
\end{ques}

In the case of $n=2$, without assuming the flatness of the normal bundle, they \cite{CMV2019} used rather different approach, namely, the Gauss-Bonnet theorem, to prove that a free boundary compact orientable surface $\Sigma^2$ immersed (resp. minimally immersed) in $\mathbb{B}^{2+k}$ with $|\Phi|^2\leq2$ (resp. $|A|^2\leq4$) is topologically a disk (resp. a flat equatorial disk).

In this paper, via some simple observations we remove the assumptions on the orientability and on the flatness of the normal bundle, answering the first question above\footnote{Very recently, Onti and Vlachos \cite{OV21} also removed the flatness assumption in their study of homology vanishing theorems for closed submanifolds by showing a generalization of the pointwise inequality (\ref{ineq-Bochner}) of Lemma \ref{lem-Bochnerineq} below for the Bochner operator. Cui and Sun \cite{CS2019} also proved Lemma \ref{lem-Bochnerineq} and used it to give a lower bound for eigenvalues of the Hodge-Laplacian on closed submanifolds.}. Moreover, we improve the pinching conditions which include the equality case of Theorem \ref{thmrefer1}, sharpen the pinching constant in the first case of Theorem \ref{thmrefer2} and include an upper bound involving the mean curvature. Besides, we also establish the vanishing theorem with the pinching condition for $|A|^2$, which can not be derived directly from that of $|\Phi|^2$ of Theorem \ref{thmrefer1}. In particular, our approach is applicable to both of the cases when $n=2$ or $n\geq3$. This is based on a new Hardy type inequality (see (\ref{Hardy1}) in Lemma \ref{Hardy-ineq}) that we establish for $n\geq2$ dimensional free boundary submanifold $M^n$ in $\mathbb{B}^{n+k}$, instead of that for $n\geq3$ of Batista-Mirandola-Vit\'{o}rio \cite{BMV17} quoted by  \cite[Lemma 2.3]{CMV2019}. This new Hardy type inequality is the key tool for improving the pinching conditions.

\begin{thm}\label{ThmA}
Let $M^n$ be a free boundary compact submanifold immersed in $\mathbb B^{n+k}$ with $n \ge 2$. If $|\Phi|^2\leq \cfrac{np}{n-p}$, for some positive integer $p \leq \lfloor \frac{n}{2} \rfloor$, then  $H^p(M; \mathbb{R}) = H^{ n-p}(M; \mathbb{R}) = 0$. In particular, if $|\Phi|^2\leq \frac{n}{n-1}$, then all cohomology groups $H^q(M; \mathbb{R})$, with $q = 1, \dots , n-1$, vanish and $M$ has only one boundary component.
\end{thm}

The pinching condition above can be easily replaced by $|A|^2\leq \cfrac{np}{n-p}$, due to $|\Phi|^2=|A|^2-n|H|^2$. With the help of the new Hardy type inequality, we are able to improve this pinching constant.

\begin{thm}\label{ThmB}
Let $M^n$ be a free boundary compact submanifold immersed in $\mathbb B^{n+k}$ with $n \ge 2$. Given a positive integer $p \leq \lfloor \frac{n}{2} \rfloor$, let $\alpha=p$ if $\alpha\leq \lfloor \frac{n}{2} \rfloor$, $\alpha=n-p$ if $\alpha> \lfloor \frac{n}{2} \rfloor$, we have the following assertions
\begin{itemize}
\item[$(1)$] For $n=2p$, $H^p(M; \mathbb{R})=0$ if $|A|^2\leq\frac{n(n+4)}{n+2}$.
\item[$(2)$] For $n>2p$, $H^{\alpha}(M; \mathbb{R})=0$  if either one of the following satisfies
\begin{itemize}
\item[$(2.1)$] $|A|^2\leq\max\limits_{s_1\leq s\leq1}s(2-s)g(s)$, in case of $s_1<1$, i.e., $\alpha=p\leq \lfloor \frac{n}{2} \rfloor$, or $\alpha=p+1=n-p$.
\item[$(2.2)$] $|A|^2\leq g(s_1)$, in case of $s_1>1$, i.e., $\alpha=n-p > \lfloor \frac{n}{2} \rfloor$ and $\alpha>p+1$.
\end{itemize}
Here $s_1$ is a constant depending on $\alpha$, and $g(s)$ is a function defined as
 $$s_1:=\frac{2\alpha}{n\beta},\quad \beta:=\frac{n-p+1}{n-p},\quad g(s):=\frac{2n}{1+\sqrt{1-\frac{4p(n-p)}{\alpha^2}(\frac{2\alpha}{n}-s)s}}.$$
\end{itemize}
\end{thm}
In the case $(1)$, the pinching constant $\frac{n(n+4)}{n+2}$ for $|A|^2$ is better than $\cfrac{np}{n-p}=n$. In particular, for $n=2$, we have shown that $M^2$ is topologically a disk if $|A|^2\leq3$ (which can also be  derived from the Gauss-Bonnet theorem as in \cite{CMV2019}).
In the case $(2.2)$, one can easily check that the pinching constant $g(s_1)>n>\cfrac{np}{n-p}$. In the case $(2.1)$, the pinching constant $\max\limits_{s_1\leq s\leq1}s(2-s)g(s)\geq \frac{2\alpha}{n}(2-\frac{2\alpha}{n})g(\frac{2\alpha}{n})>\cfrac{np}{n-p}$, which cannot be evaluated at a fixed $s$ because it varies with $n$ and $p$.

In order to give the pinching condition with the mean curvature involved, we introduce a polynomial as follows. For two parameters $\gamma\geq0$ and $0< s\leq1$, we define
\begin{equation}\label{def-Fp}
F_p(x,|H|,s,\gamma):=x^2+\frac{n(n-2p)|H|}{\sqrt{np(n-p)}}x-n|H|^2-\frac{n}{p(n-p)}\Big(s(2-s)-|H|^2\Big)\gamma.
\end{equation}
Denote by $R_p(|H|,s,\gamma)$ the positive root of $F_p(x,|H|,s,\gamma)=0$, i.e.,
\begin{equation}\label{def-Rpsr}
R_p(|H|,s,\gamma)=-\frac{n(n-2p)|H|}{2\sqrt{np(n-p)}}+\frac{n}{2\sqrt{np(n-p)}}\sqrt{n^2|H|^2+4\gamma\Big(s(2-s)-|H|^2\Big)}.
\end{equation}
Let $R_p(|H|)$ be the maximum of $R_p(|H|,s,\gamma)$ for $0< s\leq1$ and $0\leq\gamma\leq\min\{\frac{n\alpha}{2s},\frac{n^2\beta}{4}\}$, where $\alpha,\beta$ are defined in Theorem \ref{ThmB}.

\begin{thm}\label{ThmC}
With the same notation as in Theorem $\ref{ThmB}$, we have $H^{\alpha}(M; \mathbb{R})=0$ if $|\Phi|\leq R_p(|H|)$. If $|H|$ is constant, then we have
\begin{itemize}
\item[$(1)$] For $|H|\geq1$, $R_p(|H|)=R_p(|H|,s,0)=\sqrt{\frac{np}{n-p}}|H|$.
\item[$(2)$] For $|H|<1$, two subcases as in Theorem $\ref{ThmB}$ hold for $n\geq 2p$, namely,
\begin{itemize}
\item[$(2.1)$] In case of $s_1<1$, i.e., $\alpha=p\leq \lfloor \frac{n}{2} \rfloor$, or $\alpha=p+1=n-p$,
$$R_p(|H|)=\left\{\begin{array}{ll}R_p(|H|,s_1,\frac{n^2\beta}{4})& \textit{for~~} |H|<s_1,\\ R_p(|H|,|H|,\frac{n\alpha}{2|H|})& \textit{for~~} |H|\geq s_1.   \end{array}\right.$$
\item[$(2.2)$] In case of $s_1>1$, i.e., $\alpha=n-p > \lfloor \frac{n}{2} \rfloor$ and $\alpha>p+1$,
$$R_p(|H|)=R_p(|H|,1,\frac{n^2\beta}{4}).$$
\end{itemize}
\end{itemize}
In particular, if $M^n$ is minimal, i.e., $|H|\equiv0$, then the pinching condition can be rewritten as
$$|A|^2\leq\left\{\begin{array}{lll}\frac{n((n-p)^2+n)}{(n-p)(n-p+1)}& \textit{for~~} \alpha=p\leq \lfloor \frac{n}{2} \rfloor,\\ \frac{n(p(n-p)+n)}{p(n-p+1)}& \textit{for~~} \alpha=p+1=n-p,\\ \frac{n^3(n-p+1)}{4p(n-p)^2} & \textit{for~~} \alpha=n-p > \lfloor \frac{n}{2} \rfloor \textit{~~and~~} \alpha>p+1. \end{array}\right.$$
\end{thm}
The pinching constant here is better than that in the first case of Theorem \ref{thmrefer2}. As for the best pinching constants, though there certainly exist, one cannot expect them as good as those for closed minimal submanifolds in spheres. This is because that $|A|^2$ will no longer be identically equal to the pinching constants, and even in dimension $2$, as a candidate of free boundary minimal surfaces achieving the first gap for $|A|^2$, the critical catenoid has a complicated maximum of $|A|^2$.

\section{Preliminaries}
Let $M^n$ be an isometric submanifold immersed in a Riemannian manifold $N^{n+k}$ with codimension $k\ge 1$.
Let $\overline{\nabla},\nabla$ be the connections and $\overline{R}, R$ be the curvature tensors on $N$ and $M$, respectively. The second fundamental form is $A(X,Y):=\overline{\nabla}_XY-\nabla_XY$ where $X,Y$ are local vector fields on $M$. The mean curvature vector field is $H:=\frac{1}{n}\mathrm{Tr}A$, and the traceless second fundamental form is $\Phi:=A-H\mathrm{Id}$. The shape operator with respect to a normal vector $v \in T^{\bot}M$ is a self-adjoint operator $S_v:TM\rightarrow TM$ defined by $\langle S_v(X),Y\rangle=\langle A(X,Y), v\rangle$.

 The Hodge-Laplacian acting on $p$-forms of $M^n$ is defined by
$$\Delta=\mathrm d\delta+\delta d:\Omega^p(M)\rightarrow \Omega^p(M),\quad 0\leq p \leq n,$$
where $\mathrm d $ and $\delta $ are the differential and co-differential operators, respectively. The well-known Bochner formula is (cf. \cite{wuhongxi})
\begin{equation}\label{Bochner}
\Delta \omega=\nabla^*\nabla\omega+\mathcal B^{[p]}\omega.
\end{equation}
Here $\nabla ^*\nabla$ is the connection Laplacian,
and $\mathcal B ^{[p]}:\Omega^p(M)\rightarrow \Omega^p(M) $ is called the
$\mathit {Bochner\ operator}$, which can be expressed as
$$
\mathcal B^{[p]}\omega=\sum_{i,j=1}^n \theta^i\wedge i_{e_j}R(e_i,e_j)\omega,
$$
where $(e_1,\dots,e_n)$ is a local orthonormal frame and $(\theta^1,\dots,\theta^n)$ is the dual coframe. By the Gauss equation, the Bochner operator splits as (cf. \cite{Savo 2014}) $$
\mathcal B^{[p]}=\mathcal B_{\mathrm{res}}^{[p]}+\mathcal B_{\mathrm{ext}}^{[p]},
$$
where the restriction term $\mathcal B_{\mathrm{res}}^{[p]}\omega=\sum_{i,j=1}^n \theta^i\wedge i_{e_j}\overline{R}(e_i,e_j)\omega$, and the extrinsic term $\mathcal B_{\mathrm{ext}}^{[p]}$ can be expressed by the shape operators as
 $$\mathcal {B}_{\mathrm{ext}}^{[p]}=\sum_{r=1}^{k}T_{v_r}^{[p]}
=\sum_{r=1}^k\Big((\mathrm{Tr}S_{v_r})S_{v_r}^{[p]}
-S_{v_r}^{[p]}\circ S_{v_r}^{[p]}\Big).
$$
Here $(v_1,\dots,v_k)$ is a local orthonormal frame of the normal bundle $T^{\bot} M$,
$$T_v^{[p]}=(\mathrm{Tr}S_v)S_v^{[p]}-S_v^{[p]}\circ S_v^{[p]},$$
$$S_v^{[p]}\omega(X_1,\cdots,X_p)=\sum_{j=1}^p \omega(X_1,\cdots,S_v(X_j),\cdots,X_p),$$
for tangent vectors $X_1,\cdots,X_p$. The key to remove the flatness assumption of the normal bundle of Theorems \ref{thmrefer1} and \ref{thmrefer2} is the following inequality for the Bochner operator.
\begin{lem}\label{lem-Bochnerineq}
Let $M^n$ be a submanifold immersed in $N^{n+k}$ with $k\ge 1$, we have
$$\mathcal {B}_{\mathrm{ext}}^{[p]}\geq-\cfrac {p(n-p)}{n}
\left(\left|\Phi\right|^2+
\cfrac{n\left|n-2p\right|}{\sqrt{np(n-p)}}\left|H\right|\left|\Phi\right|-n\left|H\right|^2\right)\mathrm{Id}.
$$
In particular, when $N^{n+k}$ is flat, then
\begin{equation}\label{ineq-Bochner}
\langle \mathcal B^{[p]}\omega, \omega\rangle \geq-\cfrac {p(n-p)}{n}
\left(\left|\Phi\right|^2+
\cfrac{n\left|n-2p\right|}{\sqrt{np(n-p)}}\left|H\right|\left|\Phi\right|-n\left|H\right|^2\right)\left|\omega\right|^2.\end{equation}
\end{lem}
\begin{rem}
\cite{CMV2019} quoted the inequality $(\ref{ineq-Bochner})$ from \cite{lin} where the inequality was proved under the unnecessary flatness assumption of the normal bundle. As byproducts, all results in \cite{lin} hold also without this flatness assumption.
\end{rem}
\begin{proof}
For a unit normal vector $v$, we denote by $H_v:=\langle H, v\rangle$ and $\Phi_v:=\langle\Phi, v\rangle=S_v-H_v \mathrm{Id}$ the components of $H$ and $\Phi$ in the direction $v$, respectively. By choosing a principal orthonormal basis corresponding to principal curvatures $\{k_1,\cdots,k_n\}$ of $S_v$, one can show that the eigenvalues of the self-adjoint operator $T_v^{[p]}$ are given by $$\Big\{\Big(\sum_{j=1}^p k_{i_j}\Big)\Big(nH_v-\sum_{j=1}^p k_{i_j}\Big)\mid 1\leq i_1<\cdots<i_p\leq n\Big\}.$$
In \cite{Savo 2014} and \cite{lin}, it has been proven that
$$\Big(\sum_{j=1}^p k_{i_j}\Big)\Big(nH_v-\sum_{j=1}^p k_{i_j}\Big)\geq-\cfrac {p(n-p)}{n}
\Big(\left|\Phi_v\right|^2+
\cfrac{n\left|n-2p\right|}{\sqrt{np(n-p)}}\left|H_v\right|\left|\Phi_v\right|-n\left|H_v\right|^2\Big),$$
and thus
$$T_v^{[p]}\geq-\cfrac {p(n-p)}{n}
\Big(\left|\Phi_v\right|^2+
\cfrac{n\left|n-2p\right|}{\sqrt{np(n-p)}}\left|H_v\right|\left|\Phi_v\right|-n\left|H_v\right|^2\Big)\mathrm{Id}.$$

Now for a local orthonormal frame $(v_1,\dots,v_k)$  of the normal bundle, we have
$$|H|^2=\sum_{r=1}^k\left|H_{v_r}\right|^2,\quad |\Phi|^2=\sum_{r=1}^k\left|\Phi_{v_r}\right|^2.$$
Then by the Cauchy-Schwarz inequality, $\sum_{r=1}^k\left|H_{v_r}\right|\left|\Phi_{v_r}\right|\leq \left|H\right|\left|\Phi\right|$, and thus
\begin{eqnarray*}
\mathcal {B}_{\mathrm{ext}}^{[p]}=\sum\limits_{r=1}^{k}T_{v_r}^{[p]}&\geq& -\cfrac {p(n-p)}{n}\sum_{r=1}^{k}\Big(\left|\Phi_{v_r}\right|^2+
\cfrac{n\left|n-2p\right|}{\sqrt{np(n-p)}}\left|H_{v_r}\right|\left|\Phi_{v_r}\right|-n\left|H_{v_r}\right|^2\Big)\mathrm{Id}\\
&\geq&-\cfrac {p(n-p)}{n}
\Big(\left|\Phi\right|^2+
\cfrac{n\left|n-2p\right|}{\sqrt{np(n-p)}}\left|H\right|\left|\Phi\right|-n\left|H\right|^2\Big)\mathrm{Id}.
\end{eqnarray*}
When $N^{n+k}$ is flat, $\mathcal B_{\mathrm{res}}^{[p]}=0$, and thus the assertion follows.
\end{proof}

The key to improve the pinching condition is the following Hardy type inequality.
\begin{lem}\label{Hardy-ineq}
Let $M^n$ be a free boundary compact submanifold immersed in $\mathbb B^{n+k}$. For any $0< s\in\mathbb{R}$ and for any continuous function $u$ in the Sobolev spave $H_1^2(M)$, we have
\begin{eqnarray}
s(2-s)\int_M u^2-\int_M |H|^2 u^2\leq \frac{4}{n^2}\int_M|\nabla u|^2+\frac{2s}{n}\int_{\partial M}u^2, \label{Hardy1}\\
-s(2+s)\int_M u^2-\int_M |H|^2 u^2\leq \frac{4}{n^2}\int_M|\nabla u|^2-\frac{2s}{n}\int_{\partial M}u^2, \label{Hardy2}
\end{eqnarray}
where either of the equalities holds if and only if $u\equiv0$.
\end{lem}
\begin{rem}
The first inequality $(\ref{Hardy1})$ with $s=1$ reduces to that in \cite[Lemma 2.3]{CMV2019}. Obviously, $(\ref{Hardy1})$ makes sense only when $0< s\leq 1$, which will be applied to improve the pinching condition. The second inequality $(\ref{Hardy2})$ makes sense for all $s>0$ but is not applicable to the question here.
\end{rem}
\begin{proof}
Define the function $f(x):=\frac{|x|^2}{2}$ on $M^n$. Direct calculations show that $\nabla f(x)=x^{\mathrm{T}}$, where $x^{\mathrm{T}}\in T_x M$ is the tangent component of $x$ along $M$.
Moreover, the Hessian of $f(x)$ can be computed as
\begin{eqnarray*}
\mathrm{Hess}_f(X,Y)&=&\langle \nabla_X\nabla f, Y\rangle=X\langle x, Y\rangle-\langle x, \nabla_X Y\rangle\\
&=&\langle X,Y\rangle+\langle x, \overline{\nabla}_X Y-\nabla_X Y\rangle=\langle X,Y\rangle+\langle x, A(X,Y)\rangle,
\end{eqnarray*}
for tangent vector fields $X,Y$ along $x\in M$. Therefore,
$$\Delta f(x)=n+n\langle x, H\rangle.$$

Notice that on the boundary $\partial M\subset \partial \mathbb{B}^{n+k}=\mathbb{S}^{n+k-1}$ the position vector $x$ is exactly the outward normal of $\partial M$ in $M$. Then by the divergence theorem, for any function $u$ in the Sobolev spave $H_1^2(M)$, we have
\begin{eqnarray}
\int_{\partial M}u^2&=&\int_{\partial M}\langle u^2x,x\rangle=\int_M \mathrm{div}(u^2\nabla f(x))\nonumber \\
&=&\int_M \Big(\langle 2u\nabla u, x^{\mathrm{T}}\rangle +u^2\Delta f(x)\Big) \label{Hardy-eq}\\
&=&\int_M \Big(\langle 2u\nabla u, x^{\mathrm{T}}\rangle +u^2(n+n\langle x, H\rangle)\Big). \nonumber
\end{eqnarray}
For any $t>0$, by $2ab\leq ta^2+\frac{1}{t}b^2$, we have
\begin{equation}\label{Hardy-ineq11}
\left|\int_M \langle 2u\nabla u, x^{\mathrm{T}}\rangle\right|\leq \int_M 2|u| |\nabla u| |x^{\mathrm{T}}| \leq \int_M \Big(t|\nabla u|^2+\frac{1}{t}u^2|x^{\mathrm{T}}|^2\Big).
\end{equation}
Combining (\ref{Hardy-eq}) and (\ref{Hardy-ineq11}), we obtain
\begin{equation}\label{Hardy-ineq12}
\pm\left(\int_{\partial M}u^2-\int_M u^2(n+n\langle x, H\rangle)\right)\leq \int_M \Big(t|\nabla u|^2+\frac{1}{t}u^2|x^{\mathrm{T}}|^2\Big).
\end{equation}
Let $v_1:=\frac{H}{|H|}$ be the normal direction of the mean curvature at points where $H\neq0$.
Then for any $s>0$, we have
\begin{equation}\label{Hardy-ineq13}
\left|\int_M 2u^2\langle x, H\rangle\right|\leq \int_M \Big(su^2\langle x, v_1\rangle^2+\frac{1}{s}u^2|H|^2\Big).
\end{equation}
Let $t=\frac{2}{ns}$, by (\ref{Hardy-ineq12}) and (\ref{Hardy-ineq13}), we get
\begin{eqnarray}
\pm\left(\int_{\partial M}u^2-n\int_M u^2\right)&\leq& \int_M \Big(\frac{2}{ns}|\nabla u|^2+\frac{ns}{2}u^2(|x^{\mathrm{T}}|^2+\langle x,v_1\rangle^2)+\frac{n}{2s}u^2|H|^2\Big)\nonumber \\
&\leq&\int_M \Big(\frac{2}{ns}|\nabla u|^2+\frac{ns}{2}u^2|x|^2+\frac{n}{2s}u^2|H|^2\Big)\nonumber \\
&\leq&\int_M \Big(\frac{2}{ns}|\nabla u|^2+\frac{ns}{2}u^2+\frac{n}{2s}u^2|H|^2\Big),\label{Hardy-eq-cond}
\end{eqnarray}
which gives (\ref{Hardy2}) and (\ref{Hardy1}) with positive and negative signs in the left hand side, respectively. The equality condition holds because in the inequality (\ref{Hardy-eq-cond}), $|x|^2<1$ in the interior points of $M$.
\end{proof}

Now we recall harmonic forms on a manifold $M^n$ with boundary $\partial M$ (cf. \cite{Yano}).
Let $\iota : \partial M \hookrightarrow M$ be the inclusion map and $\iota^*: \Omega^*(M) \rightarrow \Omega^*(\partial M)$ be the associated pullback map of exterior forms.
For a $p$-form $\omega \in \Omega^p(M)$, the $p$-form $\iota^*\omega$ is called the \emph{tangential part}, and the $(p-1)$-form $\omega^{\bot}:=\iota^*(i_{\xi}\omega)$ is called the \emph{normal part} of $\omega$ on the boundary $\partial M$, respectively, where $\xi$ is the inward unit normal of $\partial M$ in $M$. On the boundary we have $|\omega|^2=|\iota^*\omega|^2+|\omega^{\bot}|^2$. $\omega$ is called \emph{tangential} (resp. \emph{normal}) to the boundary $\partial M$ if its normal part vanishes $\omega^{\bot}=0$ (resp. tangential part vanishes $\iota^*\omega=0$). On such $p$-forms with boundary conditions the Hodge Laplace operator $\Delta$ is self-adjoint as usually. $\omega$ is called \emph{harmonic} if $d\omega=0$ and $\delta\omega=0$.
According to Yano \cite{Yano}, harmonic $p$-forms tangential or normal to the boundary are exactly the $p$-forms in the following subspaces respectively
\begin{align*}
\mathcal H_N^p(M)&=\{\omega \in \Omega ^p(M)\mid \Delta \omega=0\ \textit{in}\ M,\ (\mathrm{d}\omega)^{\bot} =0,\ \omega^{\bot}=0\ \textit{on}\ \partial M\},\\
\mathcal H_T^p(M)&=\{\omega \in \Omega ^p(M)\mid \Delta \omega=0\ \textit{in}\ M,\ \iota^*(\delta \omega)=0,\ \iota^*(\omega)=0\ \textit{on}\ \partial M\}.
\end{align*}
Note that the subscripts $N$ and $T$ denote forms with vanishing normal and tangential part, respectively, as in that of \cite{ACS2018, CMV2019, ISS99}. $\mathcal H_N^p(M)$ (resp. $\mathcal H_T^p(M)$) corresponds to the space $\mathcal H_p^A$ (resp. $\mathcal H_p^R$) of harmonic $p$-forms with absolute (resp. relative) boundary condition in \cite{Taylor96}, where an isomorphism between $\mathcal H_p^A$ and the deRham cohomology $H^p(M; \mathbb{R})$ is proven (see also \cite{ISS99}), i.e.,
\begin{equation}\label{isom-HN-deRham}
\mathcal H_N^p(M)\cong H^p(M; \mathbb{R}).
\end{equation}
Moreover, if $M^n$ is orientable, the Hodge star operator induces an isomorphism
\begin{equation}\label{isom-HT-HN}
\mathcal H_T^{n-p}(M)\cong \mathcal H_N^p(M).
\end{equation}
An important fact is that $\dim H^{n-1}(M; \mathbb{R})\geq r-1$, where $r$ is the number of boundary components of $M$ (see \cite[Lemma 4]{ACS2018}). This implies the last conclusion of Theorems \ref{thmrefer1}, \ref{thmrefer2} and \ref{ThmA}.
The following integral version of Bochner's formula (\ref{Bochner}) (also called Weitzenb\"{o}ck's formula \cite{CMV2019, Yano} or Reilly's formula \cite{Raulot-Savo}) and refined Kato's inequality (cf. \cite{CGH00}) for harmonic forms will be useful in the proof.
\begin{lem}\label{lemBochner} \emph{(Bochner's formula)}
Let $M^n$ be a manifold with totally umbilical boundary $\partial M$ whose second fundamental form is the identity. Then we have
$$\int_M |\nabla\omega|^2+\langle \mathcal B^{[p]}\omega, \omega\rangle = -\alpha \int_{\partial M} |\omega|^2,$$
where $\alpha=p$ or $\alpha=n-p$, depending on whether $\omega\in\mathcal H_N^p(M)$ $(|\omega|^2=|\iota^*\omega|^2 \textit{~on~} \partial M )$ or $\omega\in\mathcal H_T^p(M)$ $(|\omega|^2=|\omega^{\bot}|^2 \textit{~on~}  \partial M )$, respectively.
\end{lem}
\begin{lem}\label{lemkato} \emph{(Refined Kato's inequality)}
If $\omega$ is a harmonic $p$-form on $M^n$, then
$$|\nabla\omega|^2\geq\beta |\nabla|\omega||^2,$$
where $\beta=1+\frac{1}{n-p}$ or $\beta=1+\frac{1}{p}$, depending on whether $p\leq\lfloor\frac{n}{2}\rfloor$ or $p>\lfloor\frac{n}{2}\rfloor$, respectively.
\end{lem}
Observe that one can prove the vanishing theorems without the assumption on the orientablity by only considering tangential harmonic $p$-forms in $\mathcal H_N^p(M)$ under the isomorphism (\ref{isom-HN-deRham}) for $0<p<n$. Due to the symmetry of $\alpha$ in Bochner's formula, $\beta$ in refined Kato's inequality and Lemma \ref{lem-Bochnerineq} between $p$-forms and $(n-p)$-forms, one can proceed with the proof as in the oriented case on both tangential and normal harmonic $p$-forms under the isomorphisms (\ref{isom-HN-deRham}, \ref{isom-HT-HN}) for $0<p\leq\lfloor\frac{n}{2}\rfloor$. Henceforth, we fix the range of $p$ to be $0<p\leq\lfloor\frac{n}{2}\rfloor$ and fix $\beta=1+\frac{1}{n-p}$ as in Theorems \ref{ThmA}, \ref{ThmB} and \ref{ThmC}.

\section{Proof of Theorems}
With the observation in the last section, we consider both of tangential and normal harmonic $p$-forms for $0<p\leq\lfloor\frac{n}{2}\rfloor$ as in the oriented case.
\begin{proof}[\textbf{Proof of Theorem $\mathbf{\ref{ThmA}}$}]
 Let $\omega\in\mathcal H_N^p(M)$  or $\omega\in\mathcal H_T^p(M)$ with $\alpha=p$ or $n-p$, respectively. By Lemmas \ref{lemBochner}, \ref{lemkato} and \ref{lem-Bochnerineq}, we have
\begin{eqnarray}
0&=&\alpha \int_{\partial M} |\omega|^2+\int_M |\nabla\omega|^2+\langle \mathcal B^{[p]}\omega, \omega\rangle \nonumber\\
&\geq& \alpha \int_{\partial M} |\omega|^2+\beta \int_M |\nabla|\omega||^2+\int_M\langle \mathcal B^{[p]}\omega, \omega\rangle \nonumber\\
&\geq& \alpha \int_{\partial M} |\omega|^2+\beta \int_M |\nabla|\omega||^2-\cfrac {p(n-p)}{n}\int_M F_p(\left|\Phi\right|,\left|H\right|,s,0 )
\left|\omega\right|^2, \label{ineq-main}
\end{eqnarray}
where $F_p(\left|\Phi\right|,\left|H\right|,s,0 )=\left|\Phi\right|^2+\cfrac{n(n-2p)}{\sqrt{np(n-p)}}\left|H\right|\left|\Phi\right|-n\left|H\right|^2$ as defined in (\ref{def-Fp}).
For any $t>0$, by $2ab\leq ta^2+\frac{1}{t}b^2$, we have
$$2\left|H\right|\left|\Phi\right|\leq t\left|H\right|^2+\frac{1}{t}\left|\Phi\right|^2.$$
Define functions on $t\in\mathbb{R}_+$ by
\begin{equation}\label{def-atbt}
a(t):=1+\frac{nc}{t},\quad b(t):=n(1-ct), \quad \textit{where~~} c:=\cfrac{n-2p}{2\sqrt{np(n-p)}}.
\end{equation}
 Then the inequality (\ref{ineq-main}) turns to
 \begin{eqnarray}
0&\geq& \alpha \int_{\partial M} |\omega|^2+\beta \int_M |\nabla|\omega||^2-\cfrac {p(n-p)}{n}\int_M F_p(\left|\Phi\right|,\left|H\right|,s,0 )\left|\omega\right|^2\nonumber\\
&\geq& \alpha \int_{\partial M} |\omega|^2+\beta \int_M |\nabla|\omega||^2-\cfrac {p(n-p)}{n}\int_M \left(a(t)\left|\Phi\right|^2-b(t)\left|H\right|^2\right)\left|\omega\right|^2.\label{ineq-main2}
\end{eqnarray}
Now, applying the Hardy type inequality (\ref{Hardy1}) for $0<s\leq1$ to the last term of (\ref{ineq-main2}) with $u=|\omega|$, we obtain
 \begin{eqnarray}
0&\geq& \alpha \int_{\partial M} |\omega|^2+\beta \int_M |\nabla|\omega||^2-\cfrac {p(n-p)}{n}\int_M \left(a(t)\left|\Phi\right|^2-b(t)\left|H\right|^2\right)\left|\omega\right|^2\nonumber\\
&\geq& \Big(\alpha-\cfrac {p(n-p)}{n}b(t)\frac{2s}{n}\Big) \int_{\partial M} |\omega|^2+\Big(\beta-\cfrac {p(n-p)}{n}b(t)\frac{4}{n^2}\Big) \int_M |\nabla|\omega||^2 \label{ineq-main3}\\
&&-\cfrac {p(n-p)}{n}\int_M \left(a(t)\left|\Phi\right|^2-b(t)s(2-s)\right)\left|\omega\right|^2.\nonumber
\end{eqnarray}
Consider the function $C(t):=\frac{b(t)}{a(t)}$ on $t\in\mathbb{R}_+$ with $a(t),~b(t)$ defined in (\ref{def-atbt}). For $n=2p$, $C(t)\equiv n=\frac{np}{n-p}$.
For $n>2p$, taking derivative of $C(t)$, we find that $t_0:=\sqrt{\frac{np}{n-p}}$ is the only critical point on which $C(t)$ achieves its maximum $C(t_0)=\frac{np}{n-p}$.
Now $b(t_0)=\frac{n^2}{2(n-p)}$, thus
$$\alpha-\cfrac {p(n-p)}{n}b(t_0)\frac{2s}{n}=\alpha-ps\geq \alpha-p\geq0,$$
 $$\beta-\cfrac {p(n-p)}{n}b(t_0)\frac{4}{n^2}=1+\frac{1}{n-p}-\frac{2p}{n}>0.$$
  On the other hand, $s(2-s)$ achieves its maximum $1$ at $s_0=1$. Therefore, for $t=t_0,~s=1$, if $\left|\Phi\right|^2\leq\frac{np}{n-p}$, the inequality (\ref{ineq-main3}) attains its equality. Finally, it follows from the equality condition of the Hardy type inequality (\ref{Hardy1}) that $|\omega|\equiv0$ on $M$.
  \end{proof}

\begin{proof}[\textbf{Proof of Theorem $\mathbf{\ref{ThmB}}$}]
As in the proof of Theorem \ref{ThmA} above, we have the inequality (\ref{ineq-main2}).
Since $|\Phi|^2=|A|^2-n|H|^2$ and $a(t)>0$, the inequality (\ref{ineq-main2}) implies
 \begin{eqnarray*}
0&\geq& \alpha \int_{\partial M} |\omega|^2+\beta \int_M |\nabla|\omega||^2-\cfrac {p(n-p)}{n}\int_M \left(a(t)\left|A\right|^2-(b(t)+na(t))\left|H\right|^2\right)\left|\omega\right|^2\\
&\geq& \alpha \int_{\partial M} |\omega|^2+\beta \int_M |\nabla|\omega||^2-\cfrac {p(n-p)}{n}\int_M \left(a(t)\left|A\right|^2-\widetilde{b}(t)\left|H\right|^2\right)\left|\omega\right|^2,
\end{eqnarray*}
where $\widetilde{b}(t):=b(t)+\sigma a(t)$ with $0\leq\sigma\leq n$. As before, applying the Hardy type inequality (\ref{Hardy1}) for $0<s\leq1$, we obtain
\begin{eqnarray}
0&\geq& \alpha \int_{\partial M} |\omega|^2+\beta \int_M |\nabla|\omega||^2-\cfrac {p(n-p)}{n}\int_M \left(a(t)\left|A\right|^2-\widetilde{b}(t)\left|H\right|^2\right)\left|\omega\right|^2\nonumber\\
&\geq& \Big(\alpha-\cfrac {p(n-p)}{n}\widetilde{b}(t)\frac{2s}{n}\Big) \int_{\partial M} |\omega|^2+\Big(\beta-\cfrac {p(n-p)}{n}\widetilde{b}(t)\frac{4}{n^2}\Big) \int_M |\nabla|\omega||^2 \label{ineq-main3'}\\
&&-\cfrac {p(n-p)}{n}\int_M \left(a(t)\left|A\right|^2-\widetilde{b}(t)s(2-s)\right)\left|\omega\right|^2.\nonumber
\end{eqnarray}
Consider the function $\widetilde{C}(t):=\frac{\widetilde{b}(t)}{a(t)}=C(t)+\sigma$ on $t\in\mathbb{R}_+$.
In order to get the pinching condition for $H^{\alpha}(M; \mathbb{R})=0$, we need to require that the right hand side of (\ref{ineq-main3'}) is nonnegative, which holds by requiring the following
\begin{equation}\label{pinchingcond}
\left\{\begin{array}{lll}
\widetilde{b}(t)\leq \frac{n^2}{2p(n-p)}\frac{\alpha}{s},\\
\widetilde{b}(t)\leq \frac{n^2}{2p(n-p)}\frac{n\beta}{2},\\
|A|^2\leq \widetilde{C}(t)s(2-s).
\end{array}
\right.
\end{equation}

For $n=2p$, $a(t)\equiv1$, $\widetilde{b}(t)\equiv n+\sigma$, $\widetilde{C}(t)\equiv n+\sigma$. The inequalities of (\ref{pinchingcond}) become
$$s\leq\frac{n}{n+\sigma},\quad \sigma\leq2, \quad |A|^2\leq (n+\sigma)s(2-s).$$
It is easy to see that $(n+\sigma)s(2-s)$ attains its maximum $\frac{n(n+4)}{n+2}$ when $\sigma=2$ and $s=\frac{n}{n+2}$. Therefore, the pinching condition $|A|^2\leq \frac{n(n+4)}{n+2}$ follows from the equality condition of the Hardy type inequality (\ref{Hardy1}).

For $n>2p$, we consider (\ref{pinchingcond}) according to $s_1=\frac{2\alpha}{n\beta}>1$ and $s_1<1$ respectively. Firstly, it is easy to show that $s_1\neq1$, $s_1>1$ if and only if $\alpha=n-p > \lfloor \frac{n}{2} \rfloor$ and $\alpha>p+1$, and that $s_1<1$ if and only if $\alpha=p\leq \lfloor \frac{n}{2} \rfloor$, or $\alpha=p+1=n-p$.

For $s_1>1$, noting that $s(2-s)$ attains its maximum $1$ and the second bound of (\ref{pinchingcond}) is less than the first bound  at $s=1$, we reduce (\ref{pinchingcond}) to
\begin{equation*}
\left\{\begin{array}{ll}
\widetilde{b}(t)\leq \frac{n^2}{2p(n-p)}\frac{n\beta}{2},\\
|A|^2\leq \widetilde{C}(t).
\end{array}
\right.
\end{equation*}
Since $\widetilde{C}(t)=\frac{\widetilde{b}(t)}{a(t)}$ and $\widetilde{b}(t)=b(t)+\sigma a(t)$ are decreasing on $t\geq t_0$, $\widetilde{C}(t)$ attains its maximum only if $\widetilde{b}(t)$ attains its maximum $\frac{n^2}{2p(n-p)}\frac{n\beta}{2}$, say at $t=t_1$ which depends on $0\leq\sigma\leq n$ for the moment. Meanwhile, since $a(t)$ is decreasing on $t$ and $t=t_1$ is increasing with respect to $\sigma$, $\widetilde{C}(t)$ attains its maximum at $t=t_1$ when $\sigma=n$. Namely,
$$b(t_1)+n a(t_1)=2n+nc(\frac{n}{t_1}-t_1)=\frac{n^2}{2p(n-p)}\frac{n\beta}{2}=:\frac{n\rho}{s_1},$$
where $\rho:=\frac{n\alpha}{2p(n-p)}$. It follows that
\begin{eqnarray*}
&t_1=\Big(-(\frac{\rho}{s_1}-2)+\sqrt{(\frac{\rho}{s_1}-2)^2+4nc^2}\Big)/(2c),\\
&a(t_1)=\frac{\rho}{2s_1}\Big(1+\sqrt{1-\frac{4p(n-p)}{\alpha^2}(\frac{2\alpha}{n}-s_1)s_1}\Big),\\
&\widetilde{C}(t_1)=\frac{n^2}{2p(n-p)}\frac{n\beta}{2}/a(t_1)=\frac{n\rho}{s_1}/a(t_1)=\frac{2n}{1+\sqrt{1-\frac{4p(n-p)}{\alpha^2}(\frac{2\alpha}{n}-s_1)s_1}}=g(s_1).
\end{eqnarray*}
 Thus the pinching condition for this case is $|A|^2\leq g(s_1)$.

 For $s_1<1$, we firstly observe that if $s\leq s_1$, the second bound of (\ref{pinchingcond}) is less than or equal to the first bound, and $s(2-s)$ attains its maximum at $s=s_1$. Thus it follows from the same arguments as above that the maximum of $\widetilde{C}(t)s(2-s)$ for $s\leq s_1$ is attained at $t=t_1$ and $s=s_1$, i.e.,
 $$\max_{s\leq s_1}\widetilde{C}(t)s(2-s)=s_1(2-s_1)g(s_1).$$
 On the other hand, if $s> s_1$, the second bound of (\ref{pinchingcond}) is bigger than the first bound. Thus (\ref{pinchingcond}) is reduced to
 \begin{equation*}
\left\{\begin{array}{ll}
\widetilde{b}(t)\leq \frac{n^2}{2p(n-p)}\frac{\alpha}{s},\\
|A|^2\leq \widetilde{C}(t)s(2-s).
\end{array}
\right.
\end{equation*}
As before, $\widetilde{C}(t)$ attains its maximum only if $\widetilde{b}(t)$ attains its maximum $\frac{n^2}{2p(n-p)}\frac{\alpha}{s}$ when $\sigma=n$.
Namely,
$$b(t)+n a(t)=2n+nc(\frac{n}{t}-t)=\frac{n^2}{2p(n-p)}\frac{\alpha}{s}=\frac{n\rho}{s},$$
which implies
\begin{eqnarray*}
&t=\Big(-(\frac{\rho}{s}-2)+\sqrt{(\frac{\rho}{s}-2)^2+4nc^2}\Big)/(2c),\\
&a(t)=\frac{\rho}{2s}\Big(1+\sqrt{1-\frac{4p(n-p)}{\alpha^2}(\frac{2\alpha}{n}-s)s}\Big),\\
&\widetilde{C}(t)=\frac{n^2}{2p(n-p)}\frac{\alpha}{s}/a(t)=\frac{n\rho}{s}/a(t)=\frac{2n}{1+\sqrt{1-\frac{4p(n-p)}{\alpha^2}(\frac{2\alpha}{n}-s)s}}=g(s).
\end{eqnarray*}
Hence the pinching condition for this case is $|A|^2\leq \max\limits_{s_1\leq s\leq1}s(2-s)g(s)$.
\end{proof}

\begin{proof}[\textbf{Proof of Theorem $\mathbf{\ref{ThmC}}$}]
As in the proof of Theorem \ref{ThmA}, we have the inequality (\ref{ineq-main}).
For $\gamma\geq0$, it follows from  (\ref{ineq-main}) and the Hardy type inequality (\ref{Hardy1}) for $0<s\leq1$ that
\begin{eqnarray*}
0&\geq& \alpha \int_{\partial M} |\omega|^2+\beta \int_M |\nabla|\omega||^2-\cfrac {p(n-p)}{n}\int_M F_p(\left|\Phi\right|,\left|H\right|,s,0 ) \\
&\geq& (\alpha-\frac{2s\gamma}{n})\int_{\partial M} |\omega|^2+(\beta- \frac{4\gamma}{n^2})\int_M |\nabla|\omega||^2-\cfrac {p(n-p)}{n}\int_M F_p(\left|\Phi\right|,\left|H\right|,s,\gamma )\left|\omega\right|^2,
\end{eqnarray*}
where $F_p(\left|\Phi\right|,\left|H\right|,s,\gamma )=F_p(\left|\Phi\right|,\left|H\right|,s,0 )-\frac{n}{p(n-p)}\Big(s(2-s)-|H|^2\Big)\gamma$ as defined in (\ref{def-Fp}).
Thus, in order to get the pinching condition for $H^{\alpha}(M; \mathbb{R})=0$, by the equality condition of the Hardy type inequality (\ref{Hardy1}) we only need to require
\begin{equation}\label{pinch-H}
\gamma\leq\frac{n\alpha}{2s}, \quad \gamma\leq\frac{n^2\beta}{4}, \quad F_p(\left|\Phi\right|,\left|H\right|,s,\gamma)\leq0.
\end{equation}
By definition, $$R_p(|H|)=\max\Big\{R_p(|H|,s,\gamma)\mid 0< s\leq1,~ 0\leq\gamma\leq\min\{\frac{n\alpha}{2s},\frac{n^2\beta}{4}\}\Big\},$$
where $R_p(|H|,s,\gamma)$ is the positive root of $F_p(x,|H|,s,\gamma)=0$ as expressed in (\ref{def-Rpsr}).
We remark that if $|H|$ is a non-constant function, the maximum of $R_p(|H|,s,\gamma)$ may be achieved at different values of $s$ and $\gamma$ as $|H|$ varies. So in this case, by abuse of notation we just choose fixed values $0<s=s_{\mathrm{max}}\leq1$ and $0\leq\gamma=\gamma_{\mathrm{max}}\leq \min\{\frac{n\alpha}{2s_{\mathrm{max}}}, \frac{n^2\beta}{4}\}$ such that
 $R_p(|H|)=R_p(|H|,s_{\mathrm{max}},\gamma_{\mathrm{max}})$ is as large as possible. Then the pinching condition is $|\Phi|\leq R_p(|H|)$, since now (\ref{pinch-H}) holds for $s=s_{\mathrm{max}}$ and $\gamma=\gamma_{\mathrm{max}}$. When $|H|$ is constant, we compute $s_{\mathrm{max}}$ and $\gamma_{\mathrm{max}}$ explicitly in the following.

For $|H|\geq1$, $s(2-s)-|H|^2\leq0$, thus it follows from (\ref{def-Rpsr}) that $\gamma_{\mathrm{max}}=0$ and $$R_p(|H|)=R_p(|H|,s,0)=\sqrt{\frac{np}{n-p}}|H|.$$

For $|H|<1$, two subcases as in the proof of Theorem $\ref{ThmB}$, namely, $s_1=\frac{2\alpha}{n\beta}>1$ and $s_1<1$, also occur here.
For $s_1>1$, noting that $s(2-s)$ attains its maximum $1$ and the second bound of (\ref{pinch-H}) is less than the first bound at $s=1$, we deduce immediately from (\ref{def-Rpsr}) that $s_{\mathrm{max}}=1$, $\gamma_{\mathrm{max}}=\frac{n^2\beta}{4}$ and
$$R_p(|H|)=R_p(|H|,1,\frac{n^2\beta}{4})=-\frac{n(n-2p)|H|}{2\sqrt{np(n-p)}}+\frac{n^2}{2\sqrt{np}(n-p)}\sqrt{(n-p+1)-|H|^2}.$$
For $s_1<1$, we firstly observe that if $s\leq s_1$, the second bound of (\ref{pinch-H}) is less than or equal to the first bound, and $s(2-s)$ attains its maximum at $s=s_1$. Thus
\begin{eqnarray*}
&&\max \Big\{4\gamma\Big(s(2-s)-|H|^2\Big) \mid 0< s\leq s_1,~ 0\leq\gamma\leq\min\{\frac{n\alpha}{2s}, \frac{n^2\beta}{4}\}\Big\}\\
&=&n^2\beta\Big(s_1(2-s_1)-|H|^2\Big)=2n\alpha\Big(2-(s_1+\frac{|H|^2}{s_1})\Big).
\end{eqnarray*}
On the other hand, if $s\geq s_1$, the second bound of (\ref{pinch-H}) is bigger than or equal to the first bound. Then
$$4\gamma\Big(s(2-s)-|H|^2\Big)\leq 2n\alpha\Big(2-(s+\frac{|H|^2}{s})\Big),$$
where the right hand side attains its maximum at $s=|H|$ if $|H|\geq s_1$, and at $s=s_1$ if $|H|<s_1$.
In conclusion, we have shown
$$(s_{\mathrm{max}},\gamma_{\mathrm{max}})=\left\{\begin{array}{ll}
(s_1, \frac{n^2\beta}{4}) & \textit{for~} |H|<s_1,\\
(|H|, \frac{n\alpha}{2|H|})& \textit{for~} |H|\geq s_1.
\end{array}\right.$$
Correspondingly, $R_p(|H|)=R_p(|H|,s_{\mathrm{max}},\gamma_{\mathrm{max}})$ is given by
\begin{eqnarray*}
&R_p(|H|,s_1,\frac{n^2\beta}{4})=-\frac{n(n-2p)|H|}{2\sqrt{np(n-p)}}+\frac{n^2}{2\sqrt{np(n-p)}}\sqrt{\frac{4\alpha((n-\alpha)(n-p)+n)}{n^2(n-p+1)}-\frac{|H|^2}{n-p}},\\
&R_p(|H|,|H|,\frac{n\alpha}{2|H|})=-\frac{n(n-2p)|H|}{2\sqrt{np(n-p)}}+\frac{n}{2\sqrt{np(n-p)}}\sqrt{n^2|H|^2+4n\alpha(1-|H|)}.
\end{eqnarray*}
In particular, if $M^n$ is minimal, i.e., $|H|\equiv0$, we obtain the required pinching condition
\begin{equation*}
|A|^2\leq\left\{\begin{array}{ll}
\Big(R_p(0,s_1,\frac{n^2\beta}{4})\Big)^2 & \textit{for~} s_1<1,\\
\Big(R_p(0,1,\frac{n^2\beta}{4})\Big)^2& \textit{for~} s_1>1,
\end{array}\right.
\end{equation*}
where the case of $s_1<1$ splits into two subcases when $\alpha=p\leq \lfloor \frac{n}{2} \rfloor$, or $\alpha=p+1=n-p$ in the expression $\Big(R_p(0,s_1,\frac{n^2\beta}{4})\Big)^2=\frac{n\alpha((n-\alpha)(n-p)+n)}{p(n-p)(n-p+1)}$.
\end{proof}

\begin{acknow}
The authors would like to thank Marcos P. Cavalcante for valuable discussions and for taking our attention to the paper of Onti and Vlachos.
\end{acknow}
\begin{Data}
 Data sharing not applicable to this article as no datasets were generated or analysed during the current study.
\end{Data}

\end{document}